\documentclass{amsart}

\newtheorem{theorem}{Theorem}
\newtheorem{lemma}[theorem]{Lemma}

\newcommand{\RR}{\mathbb{R}}
\newcommand{\CC}{\mathbb{C}}
\newcommand{\FF}{\mathbb{F}}

\newcommand{\ra}{\rightarrow}
\newcommand{\lra}{\longrightarrow}
\newcommand{\Rem}{\noindent {\it Remark. }}

\newcommand{\supp}{\mathord{\mathrm{supp}}}
\newcommand{\Int}{\mathord{\mathrm{Int}}}
\newcommand{\Csc}[1]{\mathord{\mathit{C}^{\infty}_{c}}(#1,\FF)}
\newcommand{\Cs}[1]{\mathord{\mathit{C}^{\infty}}(#1,\FF)}

\begin{document}

\title{On isomorphisms of algebras of smooth functions}

\author{Janez Mr\v{c}un}
\address{Department of Mathematics, University of Ljubljana,
         Jadranska 19, 1000 Ljubljana, Slovenia}
\email{janez.mrcun@fmf.uni-lj.si}
\thanks{This work was supported in part by the Slovenian Ministry of Science}
\date{}
\commby{Jonathan M. Borwein}

\begin{abstract}
We show that for any smooth Hausdorff manifolds $M$ and $N$,
which are not necessarily second-countable, paracompact or connected,
any isomorphism from the algebra of smooth
(real or complex) functions 
on $N$ to the algebra of smooth functions on $M$
is given by composition with
a unique diffeomorphism from $M$ to $N$.
An analogous result holds true for isomorphisms of algebras
of smooth functions with compact support.
\end{abstract}

\maketitle

Throughout this paper,
manifolds are assumed to be smooth and Hausdorff.
On the other hand, we do not assume manifolds
to be second-countable, paracompact or connected.
Thus manifolds may have uncountably many components,
and in particular any discrete space is a manifold of
dimension $0$.

Choose $\FF$ to be either $\RR$ or $\CC$.
For a manifold $M$ we shall denote by
$\Cs{M}$ the commutative algebra of
$\FF$-valued smooth functions on $M$.
Its subalgebra
of smooth functions with compact support
on $M$ will be denoted by $\Csc{M}$.

The aim of this paper is to prove the following
theorem, which answers the question of Alan Weinstein,
formulated at the International Euroschool on 
Poisson Geometry, Deformation Quantisation 
and Group Representations in Brussels (2003):

\begin{theorem} \label{theo1}
For any Hausdorff smooth manifolds $M$ and $N$
(which are not necessarily second-countable,
paracompact or connected),
any isomorphism of algebras of smooth functions
$T\!:\Cs{N}\ra\Cs{M}$
is given by composition with a unique diffeomorphism
$\tau\!:M\ra N$.
\end{theorem}

An analogous result, which we already proved
in \cite{Mrcun}, holds true for isomorphisms of
algebras of smooth functions with compact support:

\begin{theorem} \label{theo2}
For any Hausdorff smooth manifolds $M$ and $N$
(which are not necessarily second-countable,
paracompact or connected),
any isomorphism of algebras of smooth functions
with compact support
$\Csc{N}\ra\Csc{M}$
is given by composition with a unique diffeomorphism
$M\ra N$.
\end{theorem}

Although we have chosen to work in the smooth context,
we could as well consider manifolds of differentiability
class $\mathit{C}^{r}$, for any $r=0,1,\ldots\,$.
In that case we obtain an analogous result for the
isomorphisms of algebras of $\mathit{C}^{r}$-functions
(with compact support),
without any changes in the proof. In particular, the result
holds true for isomorphisms of algebras of continuous functions
(with compact support)
on Hausdorff topological manifolds.

The standard way to prove this kind of a statement
is to establish a one-to-one correspondence between
points of $M$ and multiplicative functionals
on the algebra of smooth functions on $M$ equipped with
a suitable topology.
However, for this one usually needs to
impose additional conditions on the spaces $M$ and $N$
(see e.g.\
\cite{Bkouche,DunfordSchwartz,GelfandKolmogoroff,Nagata,Shilov}).

In the proof of Theorem \ref{theo1}
we shall choose a different, purely algebraic approach.
We shall characterize points of a manifold by suitable
`characteristic' sequences of functions on
the manifold. We used these sequences
already in \cite{Mrcun}, where we proved
Theorem \ref{theo2}. In fact, the proof of Theorem \ref{theo1}
given below also proves Theorem \ref{theo2},
if in the argument one substitutes smooth functions
with smooth functions with compact support.

A {\em characteristic sequence of functions} on 
a manifold $M$ at a point $x\in M$
is a sequence $(f_{n})_{n}$ of $\FF$-valued smooth functions
on $M$ such that
\begin{enumerate}
\item [(i)]  $f_{n}f_{n+1}=f_{n+1}$ for any $n=1,2,\ldots\,$, and
\item [(ii)] the sequence of supports $(\supp(f_{n}))_{n}$
is a fundamental system of neighbourhoods of the point $x$ in $M$.
\end{enumerate}
Note that the equality $f_{n}f_{n+1}=f_{n+1}$ implies
that $f_{n}(x')=1$ for any $x'\in\supp(f_{n+1})$, and hence
$\supp(f_{n})$ is a neighbourhood of $\supp(f_{n+1})$.
In particular, if $(f_{n})_{n}$ is a characteristic sequence
of functions on $M$ at $x$, then $f_{n}(x)=1$ 
for any $n=1,2,\ldots\,$.
Also note that for any $x\in M$ there exists
a characteristic sequence of
functions at $x$.

\begin{lemma} \label{lem3}
Let $(f_{n})_{n}$ be a sequence of $\FF$-valued
smooth functions on a manifold $M$
satisfying $f_{n}f_{n+1}=f_{n+1}$ for any
$n=1,2,\ldots\,$.
Then $(f_{n})_{n}$ is a characteristic sequence
of functions at a point $x\in M$ if and only if
$\bigcap_{n=1}^{\infty}\supp(f_{n})=\{x\}$
and at least one of the functions in the sequence
$(f_{n})_{n}$ has compact support.
\end{lemma}

\begin{proof}
If $(\supp(f_{n}))_{n}$ is a fundamental
system of neighbourhoods of $x$, then
there exists $m\geq 1$ such that $\supp(f_{m})$
is a subset of a compact neighbourhood of $x$ and hence
itself compact. Conversely, if
$\bigcap_{n=1}^{\infty}\supp(f_{n})=\{x\}$ and
$\supp(f_{m})$ is compact for some $m\geq 1$, then
first note that $\supp(f_{n})$ is compact
for any $n\geq m$ as well.
Since for any open neighbourhood $U$ of $x$ in $M$
the intersection of the descending sequence of
compact sets
$\bigcap_{n=m}^{\infty}(\supp(f_{n})-U)$
is empty, it follows that $\supp(f_{n})-U$ is empty
for some $n\geq m$, thus $\supp(f_{n})\subset U$.
\end{proof}

The idea of the proof of Theorem \ref{theo1}
is to show that
the image of a characteristic sequence of functions
along an isomorphism of algebras is again a
characteristic sequence of functions,
and then use this to define a map between manifolds
with the required properties.
The first part of this plan
is done by the following lemma.

\begin{lemma} \label{lem4}
Let $M$ and $N$ be manifolds and let
$T\!:\Cs{N}\ra\Cs{M}$ be an isomorphism of algebras.

(i) If $(g_{n})_{n}$ is a characteristic sequence of functions
    on $N$ at a point $y\in N$, then
    $(T(g_{n}))_{n}$ is a characteristic sequence of functions
    on $M$ at a point $x\in M$.

(ii) If $(g_{n})_{n}$ and $(g'_{n})_{n}$
     are two characteristic sequences of functions
     on $N$ at the same point $y\in N$, then
     $(T(g_{n}))_{n}$ and $(T(g'_{n}))_{n}$ are characteristic
     sequences of functions on $M$ at the same point $x\in M$.
\end{lemma}

\begin{proof}
(i)
For any $n=1,2,\ldots$ write
$f_{n}=T(g_{n})$. 
Clearly we have $f_{n}f_{n+1}=f_{n+1}$,
so $f_{n}(x')=1$ for any $x'\in\supp(f_{n+1})$, and thus
$\supp(f_{n})$ is a neighbourhood of $\supp(f_{n+1})$.
Put
$$ K=\bigcap_{n=1}^{\infty}\supp(f_{n})\;.$$

First we will show that $K$ is not empty.
Suppose the opposite, that $K$ is empty.
Since each $f_{n}$ is non-zero, we can choose
a strictly increasing sequence $(i_{k})_{k}$
of natural numbers such that
$\supp(f_{i_{k}})-\supp(f_{i_{k}+1})$ is non-empty
for any $k=1,2,\ldots\,$.
This implies that
$$ U_{k}=\Int_{M}(\supp(f_{i_{k}}))-\supp(f_{i_{k}+1}) $$
is a non-empty open subset of $M$.
Since $K$ is assumed to be empty,
the disjoint sequence of open subsets
$(U_{k})_{k}$ is locally finite in $M$.
This means that we can choose a smooth function $f$
on $M$ and a sequence $(u_{k})_{k}$ of
smooth functions on $M$ such that
$\supp(u_{k})$ is a non-empty compact subset
of $U_{k}$,
$$ fu_{2k}=u_{2k} $$
and
$$ fu_{2k+1}=0 $$
for any $k=1,2,\ldots\,$. 
Note that $f_{i_{k}}u_{k}\neq 0$.
Put $g=T^{-1}(f)$ and $v_{k}=T^{-1}(u_{k})$. Since
$g_{i_{k}}v_{k}=T^{-1}(f_{i_{k}}u_{k})\neq 0$,
we can choose a point $y_{k}\in N$ such that
$$ g_{i_{k}}(y_{k})v_{k}(y_{k})\neq 0\;.$$
In particular we have $y_{k}\in\supp(g_{i_{k}})$,
so the sequence $(y_{k})_{k}$ converges to $y$
because $(g_{n})_{n}$ is a characteristic sequence
of functions at $y$. Next, we have $gv_{2k}=v_{2k}$ and
$gv_{2k+1}=0$, and since $v_{2k}(y_{2k})\neq 0$ and
$v_{2k+1}(y_{2k+1})\neq 0$, it follows
that $g(y_{2k})=1$ and $g(y_{2k+1})=0$.
As the sequence $(y_{k})_{k}$ is convergent,
this contradicts the continuity of $g$.
This proves that $K$ is not empty.

Take a point $x\in K$, and let $U$ be an
open neighbourhood of $x$ in $M$.
Choose a characteristic sequence $(\alpha_{n})_{n}$
of functions on $M$ at $x$ such that
$\supp(\alpha_{1})\subset U$.
Put $\beta_{n}=T^{-1}(\alpha_{n})$ and
$\gamma_{n}=\beta_{n}g_{n}$. Note that
$$ \gamma_{n}\gamma_{n+1}=\gamma_{n+1}\;.$$
Since $\alpha_{n}(x)=f_{n}(x)=1$, we have
$\alpha_{n}f_{n}\neq 0$ and hence
$\gamma_{n}\neq 0$.
It follows that $(\supp(\gamma_{n}))_{n}$ is a descending
sequence of non-empty sets.
Since $(g_{n})_{n}$ is a characteristic sequence,
there exists $m\geq 1$ such that $\supp(g_{n})$
is compact for any $n\geq m$. As
$\supp(\gamma_{n})\subset\supp(g_{n})$, this implies that
$\supp(\gamma_{n})$ is compact for $n\geq m$,
and therefore the intersection $\bigcap_{n}\supp(\gamma_{n})$
is a non-empty subset of $\bigcap_{n}\supp(g_{n})=\{y\}$,
so 
$$ \bigcap_{n}\supp(\gamma_{n})=\{y\}\;. $$
It now follows from Lemma \ref{lem3}
that $(\gamma_{n})_{n}$ is a characteristic
sequence of functions at $y$. In particular,
the support of $\gamma_{2}$ is a neighbourhood of $y$.
We can therefore choose $j\geq 2$ such that
$\supp(g_{j})$ is a subset of $\supp(\gamma_{2})$,
which implies that
$$ \gamma_{1}g_{j}=g_{j} $$
because $\gamma_{1}$ equals $1$ on $\supp(\gamma_{2})$
and therefore also on $\supp(g_{j})$.
This equation implies
$$ \alpha_{1}f_{j}=f_{j}\;,$$
since $T(\gamma_{1})=\alpha_{1}f_{1}$ and $f_{1}f_{j}=f_{j}$.
It follows that $\supp(f_{j})\subset\supp(\alpha_{1})\subset U$.
This shows that $(f_{n})_{n}$ is a characteristic sequence
of functions at $x$. In particular, we have $K=\{x\}$.

(ii)
For any $n=1,2,\ldots\,$ write
$f_{n}=T(g_{n})$ and $f'_{n}=T(g'_{n})$.
Part (i) implies that $(f_{n})_{n}$ is
a characteristic sequence of functions
on $M$ at a point $x\in M$, and that 
$(f'_{n})_{n}$ is
a characteristic sequence of functions
on $M$ at a point $x'\in M$. We claim that $x=x'$.
Indeed, if $x\neq x'$, we can choose $m\geq 1$ such that
$f_{m}f'_{m}=0$. This implies that $g_{m}g'_{m}=0$,
a contradiction.
\end{proof}

\begin{proof}[Proof of Theorem \ref{theo1}]
For any $x\in M$, we choose a characteristic
sequence $(f_{n})_{n}$ of functions on $M$ at $x$.
By Lemma \ref{lem4} (i) we know that
$(T^{-1}(f_{n}))_{n}$ is a characteristic sequence
of functions on $N$ at a point of $N$; denote this point by
$\tau(x)$. By Lemma \ref{lem4} (ii), this definition of $\tau(x)$
is independent of the choice of $(f_{n})_{n}$, thus we get a map
$$ \tau\!:M\lra N \;.$$

We shall now show that $T$ is given by the
composition with $\tau$. 
Take any $x\in M$, and choose a characteristic
sequence $(f_{n})_{n}$ of functions on $M$ at $x$.

First observe that if $p\in\Cs{N}$ is such that
$T(p)(x)=0$, then $p(\tau(x))=0$.
Indeed, if $p(\tau(x))\neq 0$, then $p$ has no zeros
on an open neighbourhood $V$ of $\tau(x)$ in $N$.
We can then choose
$n$ so large that $\supp(T^{-1}(f_{n}))$ is a compact
subset of $V$,
and we can define a smooth function
$q\in\Cs{N}$ with $\supp(q)\subset V$ by
$$ q(y)=\frac{T^{-1}(f_{n})(y)}{p(y)} $$
for any $y\in V$. Now we have
$$ 1=f_{n}(x)=T(T^{-1}(f_{n}))(x)=T(pq)(x)
   = T(p)(x) T(q)(x) \;,$$
so $T(p)(x)\neq 0$.

For any $g\in\Cs{N}$ we have
$$ T(g-T(g)(x) 1) (x) = T(g)(x)-T(g)(x) =0\;.$$
By the argument above this yields
$$ (g-T(g)(x)1)(\tau(x))=0\;,$$
and therefore
$$ T(g)(x)=g(\tau(x))\;.$$
This shows that $T$ is indeed given by composition
with the map $\tau$.

The uniqueness of such a map $\tau$ is clear.
Analogously, the isomorphism of algebras
$T^{-1}$ is also given by composition with
a unique map $N\ra M$,
which is the inverse of $\tau$ by uniqueness.
Since the smoothness of functions
is preserved under the composition with $\tau$
and with its inverse as well, it follows that
$\tau$ and its inverse are smooth maps.
\end{proof}

\Rem
The argument above also proves Theorem \ref{theo2},
if we use only functions with compact support.
In particular, in the proof of Theorem \ref{theo1}
we can replace the constant function $1$ on $N$
with $T^{-1}(f_{n})$.

\end{document}